\documentclass[10pt]{amsart}
\usepackage{latexsym, amsmath}

 \renewcommand{\epsilon}{\varepsilon}
 \newcommand{\newsection}[1]
 {\subsection{#1}\setcounter{theorem}{0} \setcounter{equation}{0}
 \par\noindent}

 \newtheorem{theorem}{Theorem}
 
 \newtheorem{lemma}[theorem]{Lemma}
 \newtheorem{corr}[theorem]{Corollary}
 
 \newtheorem{proposition}[theorem]{Proposition}
 \newtheorem{deff}[theorem]{Definition}
 
 \newcommand{\bth}{\begin{theorem}}
 \newcommand{\ble}{\begin{lemma}}
 \newcommand{\bcor}{\begin{corr}}
 \newcommand{\bdeff}{\begin{deff}}
 \newcommand{\bprop}{\begin{proposition}}
 \newcommand{\eth}{\end{theorem}}
 \newcommand{\ele}{\end{lemma}}
 \newcommand{\ecor}{\end{corr}}
 \newcommand{\edeff}{\end{deff}}
 
 \newcommand{\eprop}{\end{proposition}}
  \newcommand{\cd}{\, \cdot\, }

 \newcommand{\supp}{\text{supp }}
 \renewcommand{\Pi}{\varPi}
 
 \renewcommand{\Im}{\rm{Im} \,}
 \renewcommand{\epsilon}{\varepsilon}

 \newcommand{\R}{{\Bbb R}}
 
 \newcommand{\tidle}{\tilde}

\begin{document}

 \title[Global Strichartz estimates exterior to a convex obstacle]
 {Global Strichartz estimates for solutions to the wave equation exterior to a convex obstacle}
\author{Jason L. Metcalfe}
\address{School of Mathematics, Georgia Institute of Technology,
Atlanta, GA 30332} \email{metcalfe@math.gatech.edu}
\date{October 16, 2002}

\begin{abstract}
In this paper, we show that certain local Strichartz estimates for
solutions of the wave equation exterior to a convex obstacle can
be extended to estimates that are global in both space and time.
This extends the work that was done previously by H. Smith and C.
Sogge in odd spatial dimensions. In order to prove the global
estimates, we explore weighted Strichartz estimates for solutions
of the wave equation when the Cauchy data and forcing term are
compactly supported.
\end{abstract}

 \maketitle

\newsection{Introduction}
The purpose of this paper is to show that certain local Strichartz
estimates for solutions to the wave equation exterior to a
nontrapping obstacle can be extended to estimates that are global
in both space and time.  In \cite{SS1}, Smith and Sogge proved
this result for odd spatial dimensions $n\ge 3$.  Here, we extend
this result to all spatial dimensions $n\ge 3$.

If $\Omega$ is the exterior domain in $\R^n$ to a compact obstacle
and $n\ge 3$ is an even integer, we are looking at solutions to
the following wave equation

\begin{equation}\label{1.1}
\begin{cases}
\Box u(t,x)=\partial_t^2 u(t,x) - \Delta u(t,x) = F(t,x) \,, \quad
(t,x)\in \R\times\Omega\,,
\\
u(0,x)=f(x)\in \dot{H}^{\gamma}_D(\Omega)\,,
\\
\partial_tu(0,x)=g(x)\in \dot{H}^{\gamma-1}_D(\Omega)\,,
\\
u(t,x)=0 \, ,\quad x\in \partial\Omega \,.
\end{cases}
\end{equation}

\noindent Here $\Omega$ is the complement in $\R^n$ to a compact
set contained in $\{|x|\le R\}$ with $C^\infty$ boundary.
Moreover, $\Omega$ is nontrapping in the sense that there is a
$T_R$ such that no geodesic of length $T_R$ is completely
contained in $\{|x|\le R\}\cap \Omega$.  The case $\Omega=\R^n$ is
permitted.

We say that $1\le r,s\le 2\le p,q\le \infty$ and $\gamma$ are
admissible if the following two estimates hold.

\noindent {\bf Local Strichartz estimates.} {\it For
$f,g,F(t,\cdot)$ supported in $\{|x|\le R\}$, solutions to
\eqref{1.1} satisfy}
\begin{multline}\label{1.2}
\|u\|_{L^p_tL^q_x([0,1]\times\Omega)}+ \sup_{0\le t\le
1}\|u(t,\cdot)\|_{H^{\gamma}_D(\Omega)}+
\sup_{0\le t\le 1}\|\partial_t u(t,\cdot)\|_{H^{\gamma-1}_D(\Omega)}\\
\le C\,
\bigl(\,\|f\|_{H^{\gamma}_D(\Omega)}+\|g\|_{H^{\gamma-1}_D(\Omega)}
+\|F\|_{L^r_tL^s_x([0,1]\times\Omega)}\,\bigr)\,.
\end{multline}

\noindent {\bf Global Minkowski Strichartz estimates.} {\it In the
case of $\,\Omega=\R^n$, solutions to \eqref{1.1} satisfy}
\begin{multline}\label{1.3}
\|u\|_{L^p_tL^q_x(\R^{1+n})}+\sup_{t}\|u(t,\cdot)\|_{\dot{H}^{\gamma}(\R^n)}+
\sup_{t}\|\partial_t u(t,\cdot)\|_{\dot{H}^{\gamma-1}(\R^n)}
\\ \le C\,
\bigl(\,\|f\|_{\dot{H}^{\gamma}(\R^n)}+\|g\|_{\dot{H}^{\gamma-1}(\R^n)}
+\|F\|_{L^r_tL^s_x(\R^{1+n})}\,\bigr)\,.
\end{multline}

\noindent Additionally, for technical reasons we need to assume
$2>r$ and $\gamma\le (n-1)/2$. 

The global Minkowski Strichartz estimate \eqref{1.3} is a
generalization of the work of Strichartz \cite{Str1, Str2}.  The
local Strichartz estimates \eqref{1.2} for solutions to the
homogeneous ($F=0$) wave equation in a domain exterior to a convex
obstacle were established by Smith and Sogge in \cite{SS}.  In
\cite{SS1}, Smith and Sogge demonstrated that a lemma of Christ
and Kiselev \cite{CK} (see also \cite{SS1} for a proof) could be
used to establish local estimates for solutions to the
nonhomogeneous problem.

While the arguments that follow are valid in any domain exterior
to a nontrapping obstacle, it is not currently known whether the
local Strichartz estimates \eqref{1.2} hold if the obstacle is not
convex. Related eigenfunction estimates are, however, known to
fail if $\partial\Omega$ has a point of convexity.

We note here that $p,q,r,s,\gamma$ are admissible in the above
sense if the obstacle is convex, $n\ge 3$,
$$q,s'<\frac{2(n-1)}{n-3};\quad
\frac{1}{p}+\frac{n}{q}=\frac{n}{2}-\gamma=\frac{1}{r}+\frac{n}{s}-2$$
$$\frac{1}{p}=\left(\frac{n-1}{2}\right)\left(\frac{1}{2}-\frac{1}{q}\right);\quad
\frac{1}{r'}=\left(\frac{n-1}{2}\right)\left(\frac{1}{2}-\frac{1}{s'}\right)$$
where $r', s'$ represent the conjugate exponents to $r,s$
respectively.  In particular, notice that we have admissibility in
the conformal case
$$p,q=\frac{2(n+1)}{n-1};\quad r,s=\frac{2(n+1)}{n+3};\quad
\gamma=\frac{1}{2}.$$ Additionally, we note that it is well-known
(see, e.g., \cite{KT}) that in the homogeneous case ($F=0$) the
Global Minkowski Strichartz estimate \eqref{1.3} holds if and only
if $n\ge 2$, $2\le p \le \infty$, $2\le q < \infty$,
$\gamma=\frac{n}{2}-\frac{n}{q}-\frac{1}{p}$, and
\begin{equation}\label{1.4}
\frac{2}{p}\le \frac{n-1}{2}\left(1-\frac{2}{q}\right)
\end{equation}
Thus, \eqref{1.4} provides a necessary condition for
admissibility.

The main result of this paper states that for such a set of
indices a similar global estimate holds for solutions to the wave
equation in the exterior domain.
\begin{theorem}\label{theorem1.1}
Suppose $n\ge 3$.  If $p,q,r,s,\gamma$ are admissible and $u$ is a solution to the
Cauchy problem \eqref{1.1}, then
$$
\|u\|_{L^p_tL^q_x(\R\times\Omega)}\le C\,
\bigl(\,\|f\|_{\dot{H}_D^{\gamma}(\Omega)}
+\|g\|_{\dot{H}_D^{\gamma-1}(\Omega)}
+\|F\|_{L^r_tL^s_x(\R\times\Omega)}\,\bigr)\,.
$$
\end{theorem}

Throughout the sequel, we will focus on the case $n\ge 4$.  The
techniques herein can be modified to handle the $n=3$ case, but since
this case was previously handled by Smith-Sogge \cite{SS1} and since
this complicates the exposition, we choose not to provide these details.

The key differences between the general case and the odd dimensional case
are the lack of strong Huygens' principle and the fact that the
local energy no longer decays exponentially. Local energy decay
and the homogeneous Sobolev spaces $\dot{H}^{\gamma}_D(\Omega)$
will be discussed in more detail in the next section.

This paper is organized as follows.  In the next section, we will
discuss the homogeneous Sobolev spaces and the local decay of energy.
In the third section, we will establish our main estimates in
Minkowski space.  These include certain weighted Strichartz estimates
for the homogeneous wave equation with compactly supported data in
free space.  In the fourth section, we give an alternate proof of a
mixed norm estimate of Smith and Sogge \cite{SS1} which is valid in
all spatial dimensions.  Finally, in Section 5, we prove Theorem
\eqref{theorem1.1}.  

At the final stage of preparation, we learned that N. Burq
\cite{Burq} has independently obtained the results from this paper
using a slightly different method.

It is a pleasure to thank C. Sogge for his guidance and patience
during this study.  The author would also like to thank the referee
for several helpful suggestions.

\newsection{Energy Estimates}

We begin here with a few notes on the homogeneous Sobolev spaces
$\dot{H}_D^{\gamma}(\Omega)$.  The homogeneous Sobolev norms,
$\dot{H}^\gamma(\R^n)$, are given by
$$\|f\|_{\dot{H}^\gamma(\R^n)}=\|(\sqrt{-\Delta})^{\gamma}
f\|_{L^2(\R^n)}.$$
For functions supported on a fixed compact set, the homogeneous
Sobolev space $\dot{H}^\gamma(\R^n)$ are comparable to the
inhomogeneous Sobolev space $H^\gamma(\R^n)$.  Moreover, for
$\gamma<\frac{n}{2}$, the homogeneous Sobolev spaces
$\dot{H}^\gamma(\R^n)$ are preserved under multiplication by smooth
cutoff functions.

Fixing a smooth cutoff function
$\beta\in C_c^\infty$ such that $\beta(x)\equiv 1$ for $|x|\le R$,
for $|\gamma|<n/2$, we are able to define
$$\|f\|_{\dot{H}_D^\gamma(\Omega)}=\|\beta
f\|_{\dot{H}_D^{\gamma}(\tilde{\Omega})}+\|(1-\beta)f\|_{\dot{H}^\gamma(\R^n)}$$
where $\tilde{\Omega}$ is a compact manifold with boundary
containing $B_R=\Omega \cap \{|x|\le R\}$.  In particular, notice
that for functions (or distributions) supported in $\{|x|\le R\}$,
we have
$\|f\|_{\dot{H}_D^\gamma(\Omega)}=\|f\|_{\dot{H}_D^{\gamma}(\tilde{\Omega})}$.

Functions $f\in \dot{H}^\gamma_D(\tilde{\Omega})$ satisfy the
Dirichlet condition $f|_{\partial \tilde{\Omega}}=0$ (when this
makes sense). 
With the Dirichlet condition fixed, we may define the spaces
$\dot{H}_D^\gamma(\tilde{\Omega})$ in terms of eigenfunctions of
$\Delta$.  Since $\tilde{\Omega}$ is compact, we have an
orthonormal basis of $L^2(\tilde{\Omega})$, $\{u_j\}\subset
H_D^1(M)\cap C^{\infty}(M)$ with $\Delta u_j=-\lambda_j u_j$ where
$0<\lambda_j \nearrow \infty$.  Thus, for $\gamma\ge 0$, it is
natural to define
$$\dot{H}_D^\gamma(\tilde{\Omega})=\left\{v\in L^2(\tilde{\Omega})
: \sum_{j\ge 0} |\hat{v}(j)|^2 \lambda^\gamma_j <\infty\right\}$$
where $\hat{v}(j)=(v,u_j)$.  The
$\dot{H}_D^\gamma(\tilde{\Omega})$ norm is given by
$$\|v\|_{\dot{H}^\gamma_D(\tilde{\Omega})}^2 =
\sum_{j}|\hat{v}(j)|^2\lambda_j^\gamma.$$ Defining
$\dot{H}_D^\gamma(\tilde{\Omega})$ for $\gamma<0$ in terms of
duality, it is not difficult to see that the above
characterization for the norm also holds for negative $\gamma$.
Additionally, we mention that
$$\|v\|^2_{\dot{H}^1_D(\tilde{\Omega})}=\|v'\|^2_{L^2(\tilde{\Omega})}$$
and for $r<s$,
$$\|v\|^2_{\dot{H}^r_D(\tilde{\Omega})}\le
C\|v\|^2_{\dot{H}^s_D(\tilde{\Omega})}.$$ See, e.g., \cite{T2} for
further details.  Notice that by defining
$\dot{H}^\gamma_D(\tilde{\Omega})$ in this way, it builds in the
necessary compatibility conditions on the data.  For example, when
$\gamma \ge 2$, we must have
that $\Delta^j f
|_{\partial\tilde{\Omega}}=0$ for $2j\le \gamma$. 

At this point, we may define $H^\gamma_D(\Omega)$ similarly.  That is,
let
$$\|f\|_{H^\gamma_D(\Omega)}=\|\beta f\|_{H^\gamma_D(\tilde{\Omega})}
+\|(1-\beta)f\|_{H^\gamma(\R^n)}.$$
Note that since $\lambda_0>0$, we have
$\dot{H}^\gamma_D(\tilde{\Omega})=H^\gamma_D(\tilde{\Omega})$.  Also,
notice that for functions $u$ supported in $|x|<R$,
$\|u\|_{H^\gamma_D(\Omega)}\approx \|u\|_{H^\gamma_D(\tidle{\Omega})}$.

One of the key results that will allow us to establish the global
estimates from the local estimates and the global Minkowski
estimates is local energy decay.  It is this result that requires
the nontrapping assumption on the obstacle.  In odd dimensions, we
are able to get exponential energy decay:  see Taylor \cite{T},
Lax-Philips \cite{LP}, Vainberg \cite{V}, Morawetz-Ralston-Strauss
\cite{MRS}, Strauss \cite{Stra}, and Morawetz \cite{Mor1, Mor2}.
In even spatial dimensions, the decay is significantly less. The
version that we will use in this paper is

\noindent{\bf Local energy decay.} {\it For $n\ge 4$ even, data
$f,g$ supported in $\{|x|\le R\}$, $0\le\gamma$, and
$\beta(x)$ smooth, supported in $\{|x|\le R\}\,,$ there exist
$C<\infty$ such that for solutions to \eqref{1.1} where $F=0$ the
following holds
\begin{equation}\label{2.1}
\|\beta u(t,\cdot)\|_{H_D^{\gamma}(\Omega)}
+\|\beta\partial_t u(t,\cdot)\|_{H_D^{\gamma-1}(\Omega)} \le
C\, |t|^{-n/2}
\Bigl(\,\|f\|_{H_D^{\gamma}(\Omega)}+\|g\|_{H_D^{\gamma-1}(\Omega)}\,\Bigr)\,.
\end{equation}}

This is a generalized version of the results of Melrose \cite{Mel}.
Before showing how
we can derive this generalized version of local energy decay, we
would like to mention here the related works of Morawetz \cite{Mor3},
Ralston \cite{R1},
 and Strauss \cite{Stra}.

Notice that for $\gamma<n/2$, the Sobolev norms in \eqref{2.1} may be
replaced with the homogeneous Sobolev norms.  

\noindent{\it Proof of Equation \eqref{2.1}}. By density, we may,
without loss of generality, assume that $f, g$ are $C^\infty$.
When $n\ge 4$ is even, Melrose \cite{Mel} 
was able to show that a solution to the
homogeneous ($F=0$) Cauchy problem \eqref{1.1} outside a
nontrapping obstacle with data $f,g$ supported in $\{|x|\le R\}$
must satisfy
\begin{equation}\label{2.2}
\int_{B_R} |\nabla u(t,x)|^2\:dx + \int_{B_R} (\partial_t
u(t,x))^2\:dx \le C t^{-n} \left(\int |\nabla f|^2\:dx + \int
|g|^2\:dx\right)\end{equation}
 where $B_R=\{|x|\le R\}\cap
\Omega$.
Since $u$ can be controlled locally by $\nabla_x u$, \eqref{2.2}
implies 
\begin{equation}\label{2.3}
\|\beta(\cd)u(t,\cd)\|_{H^1_D(\tilde{\Omega})} +\|\beta(\cd)\partial_t
u(t,\cd)\|_{L^2(\tilde{\Omega})} \le C t^{-n/2}
\Bigl(\|f\|_{H^1_D(\tilde{\Omega})}+ \|g\|_{L^2(\tilde{\Omega})}\Bigr).
\end{equation}

Since $[\Box,\partial_t]=0$ and $\partial_t$ preserves the support of
the data and the boundary condition, we have that $u_t(t,x)$ is a
solution of
$$
\begin{cases}
\Box u_t(t,x)=0,\quad (t,x)\in \R\times\Omega\\
u_t(0,x)=g(x),\\
\partial_t u_t(0,x)=\Delta f(x),\\
u(t,x)=0, \quad x\in \partial\Omega.
\end{cases}
$$
Thus, by \eqref{2.2} and the fact that $\Box u=0$, we have
\begin{equation}\begin{split}\label{2.4}
\|\beta(\cd)u_t(t,\cd)\|_{H^1_D(\tilde{\Omega})} 
+\|\beta(\cd)\Delta &u(t,\cd)\|_{L^2(\tilde{\Omega})}\\
&=\|\beta(\cd)u_t(t,\cd)\|_{H^1_D(\tilde{\Omega})} 
+\|\beta(\cd)u_{tt}(t,\cd)\|_{L^2(\tilde{\Omega})}\\
&\le C t^{-n/2}\Bigl(\|g\|_{H^1_D(\tilde{\Omega})} + \|\Delta
f\|_{L^2(\tilde{\Omega})}\Bigr)\\
&= Ct^{-n/2}\Bigl(\|f\|_{H^2_D(\tidle{\Omega})}+\|g\|_{H^1_D(\tilde{\Omega})}\Bigr).
\end{split}\end{equation}
Thus, by elliptic regularity, \eqref{2.3}, and the monotonicity in
$\gamma$ of the norms $\|\cd\|_{H^\gamma_D(\tilde{\Omega})}$, we see
that
\begin{equation}\label{2.5}
\|\beta(\cd)u(t,\cd)\|_{H^2_D(\tilde{\Omega})}+\|\beta(\cd)\partial_t
u(t,\cd)\|_{H^1_D(\tilde{\Omega})} \le
Ct^{-n/2}\Bigl(\|f\|_{H^2_D(\tilde{\Omega})} + \|g\|_{H^1_D(\tilde{\Omega})}\Bigr).
\end{equation}

If we look similarly at $u_{tt}, u_{ttt},$ etc., we see that
\begin{equation}\label{2.6}
\|\beta(\cd)u(t,\cd)\|_{H^s_D(\tilde{\Omega})} +
\|\beta(\cd)\partial_t u(t,\cd)\|_{H^{s-1}_D(\tilde{\Omega})} \le C
t^{-n/2} \Bigl(\|f\|_{H^s_D(\tilde{\Omega})} + \|g\|_{H^{s-1}_D(\tilde{\Omega})}\Bigr)
\end{equation}
for any positive integer $s$.  By complex interpolation and the
characterization of the Sobolev spaces given above, this yields
\eqref{2.1} for any $\gamma\ge 1$.

We now work to obtain \eqref{2.1} with $\gamma < 1$.  To do so, let
$\tilde{g}$ be the solution of
$$
\begin{cases}
\Delta \tilde{g}(x)=g(x), &\text{in } \{|x|\le R\}\cap \Omega\\
\tilde{g}(x)=0, &\text{on } \{|x|=R\}\cup \partial\Omega.
\end{cases}
$$
Fix a smooth cutoff function $\chi(x)$ with $\chi(x)\equiv 1$ on
$\supp g$ and $\supp \chi \subset \{|x|<R\}$.  Then,
$$\Delta(\chi \tilde{g})= g + \psi$$
where $\psi\in C^\infty_c(\Omega)$ and by elliptic regularity,
\begin{equation}\label{2.7}
\|\psi\|_{L^2(\Omega)}\le C\|g\|_{H^{-1}_D(\Omega)}.
\end{equation}
If $v$ is the solution to
$$
\begin{cases}
\Box v(t,x)=0\\
v(0,x)=\chi(x) \tilde{g}(x)\\
\partial_t v(0,x)=f(x)\\
v(t,x)=0, \quad x\in\partial\Omega,
\end{cases}
$$
since $[\Box,\partial_t]=0$, we have
$$
\begin{cases}
\Box \partial_t v(t,x)=0\\
\partial_t v(0,x)=f(x)\\
\partial_t^2 v(0,x)=\Delta (\chi(x)\tilde{g}(x))\\
\partial_t v(t,x)=0, \quad x\in\partial\Omega.
\end{cases}
$$
Thus, by \eqref{2.2} and \eqref{2.7}, we have
\begin{equation}\begin{split}\label{2.8}
\|\beta(\cd)\partial_t v(t,\cd)\|_{L^2(\tilde{\Omega})}+\|\beta(\cd)&\partial_t^2
v(t,\cd)\|_{H^{-1}_D(\tidle{\Omega})}\\
&\le C \|\beta(\cd)\Delta v(t,\cd)\|_{H^{-1}_D(\tilde{\Omega})}+\|\beta(\cd)\partial_t
v(t,\cd)\|_{L^2(\tidle{\Omega})}\\
&\le C \|\beta(\cd)v(t,\cd)\|_{H^1_D(\tilde{\Omega})}+\|\beta(\cd)\partial_t
v(t,\cd)\|_{L^2(\tilde{\Omega})}\\
&\le C t^{-n/2}
\Bigl(\|\chi(x)\tilde{g}(x)\|_{H^1_D(\Omega)}+\|f\|_{L^2(\Omega)}\Bigr)\\
&\le C t^{-n/2}\Bigl(\|f\|_{L^2(\Omega)}+\|g\|_{H^{-1}_D(\Omega)}\Bigr).
\end{split}
\end{equation}

Since $u-\partial_t v$ also solves a homogeneous wave equation with
$C^\infty_c(\Omega)$ data
$$
\begin{cases}
\Box (u-\partial_t v)=0\\
(u-\partial_t v)(0,\cd)=0\\
(\partial_t u - \partial_t^2 v)(0,\cd)=\psi\\
(u-\partial_t v)(t,x)=0, \quad \text{for } x\in\partial\Omega,
\end{cases}
$$
we have
\begin{equation}\begin{split}\label{2.9}
\|\beta(\cd)(u-\partial_t
v)(t,\cd)&\|_{L^2(\tilde{\Omega})}+\|\beta(\cd) (\partial_t u - \partial_t^2
v)(t,\cd)\|_{H^{-1}_D(\tilde{\Omega})} 
\\&\le \|\beta(\cd)(u-\partial_t
v)(t,\cd)\|_{H^1_D(\tilde{\Omega})}+\|\beta(\cd) (\partial_t u - \partial_t^2
v)(t,\cd)\|_{L^2(\tilde{\Omega})}\\
&\le C t^{-n/2} \|\psi\|_{L^2(\Omega)}\\
&\le C t^{-n/2} \|g\|_{H^{-1}_D(\Omega)}.
\end{split}\end{equation}
Combining \eqref{2.8} and \eqref{2.9}, it follows easily that
\begin{equation}\label{2.10}
\|\beta(\cd) u(t,\cd)\|_{L^2(\tilde{\Omega})}+ \|\beta(\cd)\partial_t
u(t,\cd)\|_{H^{-1}_D(\tilde{\Omega})} \le C t^{-n/2}\Bigl(\|f\|_{L^2(\Omega)}+\|g\|_{H^{-1}_D(\Omega)}\Bigr).
\end{equation}
Finally, if we interpolate with \eqref{2.3}, we see that we obtain
\eqref{2.1} for $0\le \gamma \le 1$ which completes the proof. \qed

\newsection{Weighted Minkowski Estimates}

In this section we show that weighted versions of the Minkowski
Strichartz estimates for solutions to the homogeneous wave
equation can be obtained when the initial data are compactly
supported.  Specifically, we are looking at the homogeneous free
wave equation
\begin{equation}\label{free}
\begin{cases}
\Box w(t,x)=\partial_t^2 w(t,x) - \Delta w(t,x) = 0 \,, \quad
(t,x)\in \R\times\R^n\,,
\\
w(0,x)=f(x)\in H^{\gamma}(\R^n)\,,
\\
\partial_tw(0,x)=g(x)\in H^{\gamma-1}(\R^n)\,.
\end{cases}
\end{equation}
where the Cauchy data $f,g$ are supported in
$\{x\in\R^n:|\,x\,|<R\}$.


We begin by showing that one can obtain weighted versions of the
energy inequality.  Here we need only slightly
modify the arguments of H\"{o}rmander \cite{H} (Lemma 6.3.5, p.
101) and Lax-Philips \cite{LP} (Appendix 3).

\begin{lemma}\label{lemma3.1}
Suppose that $n\geq 3$. Let $w(t,x)$ be a solution to the
homogeneous Minkowski wave equation \eqref{free} with smooth
initial data $f,g$ supported in $\{|x|\leq R\}$. Then, the
following estimate holds
$$\int (t-|x|)^2\left(\bigl|\nabla_x w(t,x)\bigr|^2 + (\partial_t w(t,x))^2\right)\:dx
\leq C_R \left(\int \bigl|\nabla f\bigr|^2 + |g|^2\:dx\right)$$
\end{lemma}

\noindent{\it Proof.} It is not difficult to check that
 $$\text{div}_x p + \partial_t q = N(w)\Box w$$
where
\begin{align*}
N(w)&=4t(x\cdot \nabla w) + 2(r^2+t^2)w_t + 2(n-1)tw\\
p&=-2t w_t^2 x - 4t(x\cdot \nabla w)\nabla w + 2 t \bigl|\nabla
w\bigr|^2 x\\
&\quad\quad\quad\quad\quad\quad - 2(r^2+t^2)w_t\nabla w - 2(n-1)t w \nabla w\\
q&=4t(x\cdot\nabla w)w_t + (r^2+t^2)\bigl(\bigl|\nabla
w\bigr|^2+w_t^2\bigr)+ 2(n-1)tww_t - (n-1)w^2.
\end{align*}
If we integrate over a cylinder $[0,T]\times \{x\in\R^n\,:\,
|x|\le \bar{R}\}$ for $\bar{R}$ sufficiently large, Huygens'
principle and the divergence theorem gives us that:
\begin{equation}\label{3.2}
\int_{t=T} q\:dx - \int_{t=0} q\: dx = 0.
\end{equation}
Here, since the initial data are compactly supported, we have
\begin{equation}\label{3.3}
\begin{split}
\int_{t=0} q\: dx &= \int_{t=0} r^2\bigl(\bigl|\nabla_x
w(0,x)\bigr|^2
+ w_t(0,x)^2\bigr) - (n-1) w(0,x)^2 \:dx\\
&\leq C_R \left(\int \bigl|\nabla f\bigr|^2 + |g|^2\:dx\right).
\end{split}
\end{equation}

Now, let us introduce the standard invariant vector fields
$$Z_0= t\partial_t + \sum_{j=1}^n x_j\partial_j, \quad
Z_{0k}=t\partial_k + x_k\partial_t, \quad
 Z_{jk}=x_k\partial_j -
x_j\partial_k$$
 for $j,k=1,2,...,n$.  Notice that
\begin{equation}\label{3.4}
\int_{t=T} q\:dx =\int_{t=T}\left( \bigl|Z_0w\bigr|^2 + \sum_{0\le
j<k\le
  n}\bigl|Z_{jk}w\bigr|^2
+2(n-1)tww_t - (n-1)w^2\right)\:dx\end{equation}

Applying Lemma 6.3.5 of H\"{o}rmander \cite{H} (p. 101), we see
that \eqref{3.2}-\eqref{3.4} yield
\begin{equation}\label{3.5}
\|Z_0w(t,\cdot)\|^2_{L^2(\R^n)}+\sum_{j<k}\|Z_{jk}w(t,\cdot)\|^2_{L^2(\R^n)}
\le C_R\left(\int \bigl|\nabla f\bigr|^2 + |g|^2\:dx\right).
\end{equation}

\noindent Thus, we see that in order to complete the proof, it
suffices to show that
$$\int (t-r)^2\left(w_t^2+\bigl|\nabla w\bigr|^2\right)\:dx\leq
\|Z_0w\|^2_{L^2(\R^n)}+ \sum_{0\leq j<k\leq
n}\|Z_{jk}w\|^2_{L^2(\R^n)}.$$

Since the Cauchy-Schwarz inequality gives us that $\bigl|\nabla
w\bigr|\geq w_r$ and since $4trw_tw_r\ge -2tr(w_t^2+w_r^2)$, we
have
\begin{align*}
\|Z_0w\|^2_{L^2(\R^n)}+ &\sum_{0\leq j<k\leq
  n}\|Z_{jk}w\|^2_{L^2(\R^n)}
=\int (t^2+r^2)\left(w_t^2+\bigl|\nabla
  w\bigr|^2\right)+4trw_tw_r\:dx\\
&=\int (t-r)^2\left(w_t^2+\bigl|\nabla
  w\bigr|^2\right)+2trw_t^2+2tr\bigl|\nabla w\bigr|^2 + 4trw_tw_r \\
&\geq \int (t-r)^2\left(w_t^2+\bigl|\nabla
  w\bigr|^2\right)
\end{align*}
as desired. \qed


Next, we look at the weighted analog of the dispersive inequality
when the initial data have compact supports.  

\begin{lemma}\label{lemma3.3}  Suppose $n\ge 2$.
Let $w$ be a solution to the homogeneous Minkowski wave equation
\eqref{free} with initial data $f,g$ supported in $\{|x|\leq R\}$.
Then, we have
\begin{equation}\label{3.6}
\|(|t|-|x|)^{(n-1)/2} w(t,x)\|_{L^{\infty}_x(\{|t|-|x|\ge 2R\})}\le
\frac{C_R}{|t|^{(n-1)/2}}\left(\|f\|_{L^2(\R^n)}+\|g\|_{L^2(\R^n)}\right),
\end{equation}
\begin{equation}\label{3.7}
\|(|t|-|x|)^{(n+1)/2} \partial_t w(t,x)\|_{L^\infty_x(\{|t|-|x|\ge
2R\})}\le
\frac{C_R}{|t|^{(n-1)/2}}\left(\|f\|_{L^2(\R^n)}+\|g\|_{L^2(\R^n)}\right).
\end{equation}
Additionally, for any $n\ge 4$,
\begin{equation}\label{3.8}
\|(|t|-|x|)^\theta w(t,x)\|_{L^2_x (\{|t|-|x|\ge 2R\})}\le
C_{R,\theta} \left(\|f\|_{L^2(\R^n)} + \|g\|_{L^2(\R^n)}\right)
\end{equation}
for any $\theta < 1$.
\end{lemma}

We note that \eqref{3.8} holds for any $\theta < 1/2$ when $n=3$.  
This is sufficient to yield the results in the
sequel for $q>2$.  Since, however, the $n=3$ case was handled by Smith and
Sogge \cite{SS1} and since this would complicate the argument, we
choose not to provide the details here.

\noindent{\it Proof.} By scaling, we may assume that $R=1$.  For
simplicity, we will demonstrate the result for $t>0$.

Begin by writing $w=w_1+w_2$, where $w_1$ is a solution of the
homogeneous Minkowski wave equation \eqref{free} with Cauchy data
$(w,w_t)|_{t=0}=(f,0)$ and $w_2$ is a solution of the Minkowski
wave equation \eqref{free} with Cauchy data
$(w,w_t)|_{t=0}=(0,g)$. It will, thus, suffice to show that the
estimate holds for $w_1$ and $w_2$ separately. Since the arguments
are the same for each piece, we will restrict our attention to
showing that the estimate holds for $w_2$, the more technical
piece.

From equations (5.43) and (5.48) of \cite{T2} (p. 222), we have that
\begin{equation}\label{3.9}
w_2(t,x)=R(t,\cd)*g
\end{equation}
where
\begin{equation}\label{3.10}
R(t,x)=\lim_{\varepsilon \searrow 0} c_n \Im (|x|^2-(t-i\varepsilon)^2)^{-(n-1)/2}.
\end{equation}
Thus, since $g$ is supported in $\{|x|<1\}$, we can apply the Schwarz
inequality to see
$$
|w(t,x)| \le C \Bigl(\int_{|y|\le 1}
 |R(t,x-y)|^2\:dy\Bigr)^{1/2} \|g\|_{L^2(\R^n)} \le C_n \sup_{|y|\le
1} |R(t,x-y)|\: \|g\|_{L^2(\R^n)}.$$
Since $|x|\le |t|-2$ and $|y|\le 1$, we have that $t^2-|x-y|^2 \approx
t^2-|x|^2$.  Thus, by \eqref{3.10}, we have
$$\sup_{|y|\le 1}|R(t,x-y)| \le C |t|^{-(n-1)/2}
(|t|-|x|)^{-(n-1)/2}$$
which completes the proof of \eqref{3.6}.  Since it is easy to see 
$$\sup_{|y|\le 1}|\partial_t R(t,x-y)| \le C
|t|^{-(n-1)/2}(|t|-|x|)^{-(n+1)/2},$$
we also get \eqref{3.7}.

For \eqref{3.8}, we again use \eqref{3.9} and \eqref{3.10} to see that
\begin{multline*}
\|(|t|-|x|)^\theta w(t,x)\|_{L^2_x(\{|t|-|x|\ge 2\})} \\\le
C \Bigl\|(t-|x|)^\theta \int \frac{1}{(t^2-|x-y|^2)^{(n-1)/2}}
g(y)\:dy\Bigr\|_{L^2(\{t-|x|\ge 2\})}.
\end{multline*}
Since $|y|\le 1$ and $t-|x|\ge 2$, we have that the right hand side is
controlled by
$$C \Bigl\|\int \frac{1}{(t+|x-y|)^{(n-1)/2}
(t-|x-y|)^{((n-1)/2)-\theta}} g(y)\:dy\Bigr\|_{L^2(\{t-|x|\ge 2\})}.$$
By Young's inequality, this is dominated by
$$C
\Bigl\|\frac{1}{(t+|x|)^{(n-1)/2}}\frac{1}{(t-|x|)^{((n-1)/2)-\theta}}\Bigr\|_{L^2(\{t-|x|\ge
1\})} \|g\|_1.$$
Since $g$ is compactly supported, by the Schwarz inequality, we have
that $\|g\|_1\le C \|g\|_2$.  We thus want to examine the $L^2$ norm
of the kernal above.  Writing this in polar coordinates, we see that
the square of this norm is bounded by
$$\int_0^{t-1}  \int_{S^{n-1}} \frac{1}{(t-\rho)^{n-1-2\theta}}\:d\sigma(\omega)\:d\rho.$$
This establishes \eqref{3.8} since this integral is bounded
independent of $t$ for $n\ge 4$ and any $\theta < 1$. \qed


From the previous three lemmas, we are able to derive a weighted
Strichartz estimate for solutions to the Minkowski wave equation
with compactly supported initial data.

\begin{theorem}\label{theorem3.4}
Suppose $n\ge 4$ and $p,q,\gamma$ are admissible.  Let $w$ be a
solution to the homogeneous Minkowski wave equation \eqref{free}
with Cauchy data $f,g$ supported in $\{|x|\le R\}$.  Then, for any
$\theta<1$, we have the following estimate:
$$\|(|t|-|x|)^\theta w(t,x)\|_{L^p_tL^q_x(\R^{1+n})}\le
C_R\bigl(\|f\|_{H^\gamma(\R^n)}+
\|g\|_{H^{\gamma-1}(\R^n)}\bigr).$$
\end{theorem}

\noindent{\it Proof.} By the Global Minkowski Strichartz estimate
\eqref{1.3} and finite propogation speed, it will suffice to show the
estimate in the case $|t|-|x|\ge 2R$.  We will, also, stick to the
case $t\ge 0$. Let $S_t=\{x: t-|x|\ge 2R\}$.

By Lemma \ref{lemma3.1} and Lemma \ref{lemma3.3}, we have:
\begin{align*}
\|(t-|x|)^\theta\partial_t w(t,x)\|_{L^2_x(S_t)}&\le C\left(\|f\|_{H^1(\R^n)}+\|g\|_{L^2(\R^n)}\right)\\
\|(t-|x|)^\theta\partial_t w(t,x)\|_{L^\infty_x(S_t)}&\le
\frac{C}{t^{(n-1)/2}}\left(\|f\|_{H^1(\R^n)}+\|g\|_{L^2(\R^n)}\right).
\end{align*}
In the second inequality, we have used the monotonicity in $\gamma$ of
$H^{\gamma}$. By Riesz-Thorin interpolation, we
have
$$
\|(t-|x|)^\theta\partial_t w(t,x)\|_{L^q_x(S_t)}\le C\:
\left(\frac{1}{t^{(n-1)/2}}\right)^{\left(1
-\frac{2}{q}\right)}\left(\|f\|_{H^1(\R^n)}+\|g\|_{L^2(\R^n)}\right).
$$

Since by \eqref{1.4}
$$p\cdot \frac{n-1}{2}\left(1-\frac{2}{q}\right) >
\frac{p}{2}\left(\frac{n-1}{2}\Bigl(1-\frac{2}{q}\Bigr) \right)
\ge
  1,$$
we see that taking the $L^p_t$ norm of both sides yields
\begin{equation}\label{3.12}
\|(t-|x|)^\theta \partial_t w(t,x)\|_{L^p_tL^q_x(\{t\ge 2R\}\times
S_t)}\le C
  \left(\|f\|_{H^1(\R^n)}+\|g\|_{L^2(\R^n)}\right)
  \end{equation}
for $n\ge 3$.

Similarly, we may interpolate between \eqref{3.6} and \eqref{3.8} to
see that
\begin{equation}\label{3.13}
\|(t-|x|)^\theta w(t,x)\|_{L^p_t L^q_x (\{t\ge 2R\}\times S_t)} \le C
\left(\|f\|_{L^2(\R^n)} + \|g\|_{L^2(\R^n)}\right)
\end{equation}
for $n\ge 4$.

If we now argue as we did in obtaining \eqref{2.10} from
\eqref{2.3}, we see that \eqref{3.12} and \eqref{3.13} yield
\begin{equation}\label{3.14}
\|(t-|x|)^\theta w(t,x)\|_{L^p_tL^q_x(\{t\ge 2R\}\times S_t)}\le C
  \left(\|f\|_{L^2(\R^n)}+\|g\|_{\dot{H}^{-1}(\R^n)}\right).
\end{equation}
The result, then, follows from the monotonicity of the Sobolev
norms.
\qed

\newsection{Mixed Estimates in Minkowski Space}

In this section, as in Smith and Sogge \cite{SS1}, we collect a couple of results that follow from the
fact that
\begin{equation}\label{4.1}
\sup_\xi\, |\xi|^{2\gamma} \left[\,
\int\bigl|\widehat\beta(\xi-\eta)\bigr|\;\delta(\tau-|\eta|\,)\,d\eta
\,\right]\le C_{n,\gamma,\beta}\,\tau^{2\gamma}
\end{equation}
if $\beta$ is a smooth function supported in $\{|x|\leq 1\}$ and
$0\leq \gamma\leq\frac{n-1}{2}$.

The first of these results is Lemma 2.2 of \cite{SS1}. 
\begin{lemma}\label{lemma4.1}\
Let $\beta$ be a smooth function supported in $\{|x|\leq 1\}$.
Suppose $0\le \gamma\le\frac{n-1}{2}$.  Then
$$\int_{-\infty}^{\infty}\bigl\|\beta(\cdot) e^{it\sqrt{-\Delta}}f(\cdot)
\bigr\|^2_{\dot{H}^{\gamma}(\R^n)}\:dt\leq
C_{n,\gamma,\beta}\,\|f\|^2_{\dot{H}^{\gamma}(\R^n)}$$
\end{lemma}

The second result is an analog of a result in \cite{SS1} that was
shown in odd spatial dimensions.  Here we require a different argument
that does not rely on sharp Huygens' principle.

\begin{lemma}\label{lemma4.2}
Let $w$ be a solution to the Cauchy problem for the Minkowski wave
equation
$$\begin{cases}
\Box w(t,x)=\partial_t^2 w(t,x) - \Delta w(t,x) = F(t,x) \,, \quad
(t,x)\in \R\times\R^n\,,
\\
w(0,x)=f(x)\,,
\\
\partial_tw(0,x)=g(x)\,.
\end{cases}$$
  Suppose that the global Minkowski Strichartz
estimate $(1.3)$ holds, that $0\le\gamma \le \frac{n-1}{2}$, and that
$r<2$. Then, for $\beta$ a smooth function supported in
$\{\,|x|\le 1\,\}\,$, we have
$$\sup_{|\alpha|\le 1}\; \int_{-\infty}^{\infty} \bigl\|\,\beta(\cdot)\,
\partial_{x}^\alpha w(t,\cdot)\,
\bigr\|^2_{\dot{H}^{\gamma-1}(\R^n)}\:dt \le
C\,\bigl(\,\|f\|_{\dot{H}^\gamma(\R^n)} +
\|g\|_{\dot{H}^{\gamma-1}(\R^n)} +
\|F\|_{L^r_tL^s_x(\R^{1+n})}\,\bigr)^2.$$
\end{lemma}

\noindent{\it Proof.} If $F=0$, the result follows from Lemma
\ref{lemma4.1}.  Thus, it will suffice to show that
$$\int_0^{\infty} \bigl\|\,\beta(\cdot)\,w(t,\cdot)\bigr\|^2_{\dot{H}^{\gamma}(\R^n)}\:dt
\le C\,\|F\|^2_{L^r_tL^s_x(\R^{1+n})}$$ when the initial data
$f,g$ are assumed to vanish.
\par
We begin by establishing that
$$TF(t,x)=\Lambda^{\gamma}\beta(\cdot)\int\frac{\sin(t-s)\Lambda}{\Lambda}F(s,\cdot)\:ds$$
is bounded from $L^r_tL^s_x(\R_+^{1+n})$ to
$L^2_tL^2_x(\R^{1+n})$. In other words, we want to show that
\begin{equation}\label{4.3}
\int
\bigl\|\beta(\cdot)\int\frac{\sin(t-s)\Lambda}{\Lambda}F(s,\cdot)
\:dx\bigr\|^2_{\dot{H}^\gamma(\R^n)}\:dt\le
C\|F\|^2_{L^r_tL^s_x(\R^{1+n})}
\end{equation}
when $F$ is assumed to vanish for $t<0$.
\par
By Strichartz estimate \eqref{1.3}, we have
\begin{equation}
\begin{split}\label{4.4}
\int\frac{|\eta|^{2\gamma}}{|\eta|^2}|\tilde{F}(|\eta|,\eta)|^2\:d\eta
&=
\int \frac{|\eta|^{2\gamma}}{|\eta|^2}\left|\int e^{-is|\eta|}\hat{F}(s,\eta)\:ds\right|^2\:d\eta\\
&\le \sup_t \int |\eta|^{2\gamma}\left|\int_0^t
\frac{e^{i(t-s)|\eta|}}{|\eta|}\hat{F}(s,\eta)\:ds\right|^2
\:d\eta\\
&\le \sup_t \left(\|w(t,\cdot)\|^2_{\dot{H}^{\gamma}(\R^n)}+\|\partial_t w(t,\cdot)\|^2_{\dot{H}^{\gamma-1}(\R^n)}\right)\\
&\le C\|F\|^2_{L^r_tL^s_x(\R^{1+n})}
\end{split} \end{equation}
where $\tilde{F}$ denotes the space-time Fourier transform of $F$.
\par
By Plancherel's theorem in $t,x$, we have
\begin{multline*}
\int
\bigl\|\beta(\cdot)\int\frac{e^{i(t-s)\Lambda}}{\Lambda}F(s,\cdot)
\:ds\bigr\|^2_{\dot{H}^\gamma(\R^n)}\:dt =\\ \int_0^\infty \int
|\xi|^{2\gamma}\left|\int\hat{\beta}(\xi-\eta)
\delta(\tau-|\eta|)\frac{1}{|\eta|}\tilde{F}(|\eta|,\eta)\:d\eta\right|^2\:d\xi\:d\tau.
\end{multline*}
By the Schwarz inequality in $\eta$, this can be bounded by
\begin{multline*}\int_0^\infty\int |\xi|^{2\gamma}\left(\int
|\hat{\beta}(\xi-\eta)|\delta(\tau-|\eta|)\:d\eta\right)\\ \times
\left(\int
|\hat{\beta}(\xi-\eta)|\delta(\tau-|\eta|)\frac{1}{|\eta|^2}|\tilde{F}(|\eta|,\eta)|^2\:d\eta
\right)\:d\xi\:d\tau.\end{multline*}

Applying \eqref{4.1}, \eqref{4.4}, and Young's inequality, we see
that
\begin{align*}
\int
\bigl\|\beta(\cdot)&\int\frac{e^{i(t-s)\Lambda}}{\Lambda}F(s,\cdot)
\:ds\bigr\|^2_{\dot{H}^\gamma(\R^n)}\:dt\\
&\le \int_0^{\infty}\int \tau^{2\gamma}\int
|\hat{\beta}(\xi-\eta)|\delta(\tau-|\eta|)\frac{1}{|\eta|^2}
|\tilde{F}(|\eta|,\eta)|^2\:d\eta\:d\xi\:d\tau\\
&=\int\int
|\hat{\beta}(\xi-\eta)|\frac{|\eta|^{2\gamma}}{|\eta|^2}
|\tilde{F}(|\eta|,\eta)|^2\:d\eta\:d\xi\\
&\le C\int\frac{|\eta|^{2\gamma}}{|\eta|^2}|\tilde{F}(|\eta|,\eta)|^2\:d\eta\\
&\le C\|F\|^2_{L^r_tL^s_x(\R^{1+n})}.
\end{align*}
By a similar argument, we can show that the same bound holds for
$$\int \bigl\|\beta(\cdot)\int\frac{e^{-i(t-s)\Lambda}}{\Lambda}F(s,\cdot)
\:ds\bigr\|^2_{\dot{H}^\gamma(\R^n)}\:dt$$ which establishes
\eqref{4.3}.  By duality, this is equivalent to having
$$T^*F:L^2_tL^2_x(\R^{1+n})\to L^{r'}_tL^{s'}_x(\R^{1+n})$$
bounded, where
$$T^*F=\int \frac{\sin(s-t)\Lambda}{\Lambda}\beta(\cdot)\Lambda^{\gamma}F(s,\cdot)\:ds.$$
\par
We wanted to show, instead, that
$$WF:L^r_tL^s_x(\R^{1+n})\to L^2_tL^2_x(\R^{1+n})$$
is bounded, where
$$WF(t,x)=\Lambda^{\gamma}\beta(\cdot)\int_0^t \frac{\sin(t-s)\Lambda}{\Lambda}F(s,\cdot)\:ds.$$
By duality, this is equivalent to showing that
$$W^{*}F:L^2_tL^2_x(\R^{1+n})\to L^{r'}_tL^{s'}_x(\R^{1+n})$$
where
$$W^{*}F(t,x)=\int_t^{\infty}\frac{\sin(s-t)\Lambda}{\Lambda}\beta(\cdot)\Lambda^{\gamma}F(s,\cdot)\:ds.$$
This, however, follows from \eqref{4.3} after an application of
the following lemma of Christ and Kiselev \cite{CK} (see also
\cite{SS1}). \qed

\begin{lemma}\label{lemma4.3}
Let $X$ and $Y$ be Banach spaces and assume that $K(t,s)$ is a
continuous function taking its values in $B(X,Y)$, the space of
bounded linear mappings from $X$ to $Y$.  Suppose that $-\infty\le
a < b \le \infty$ and $1\le p < q\le\infty$.  Set
$$Tf(t)=\int_a^b K(t,s)f(s)\:ds$$
and
$$Wf(t)=\int_a^t K(t,s)f(s)\:ds.$$
Suppose that
$$\|Tf\|_{L^q([a,b],Y)}\le C\|f\|_{L^p([a,b],X)}.$$
Then,
$$\|Wf\|_{L^q([a,b],Y)}\le C\|f\|_{L^p([a,b],X)}.$$
\end{lemma}

\newsection{Strichartz Estimates in the Exterior Domain}

By scaling, we may take $R=\frac{1}{2}$ in the sequel.  We begin
by proving a weighted version of Theorem \ref{theorem1.1} when the
data and forcing terms are compactly supported.

\begin{lemma}\label{lemma5.1}
Suppose $n\ge 4$, and suppose $u$ is a solution to the Cauchy problem \eqref{1.1} with
the forcing term $F$ replaced by $F+G$, where $F,G$ are supported
in $\{|t|\le 1\}\times \{|x|\leq 1\}$ and the initial data $f,g$
are supported in $\{|x|\leq 1\}$.  Then, for admissible
$p,q,r,s,\gamma$, there exist a positive, finite constant $C$ 
so that the following estimate holds:
\begin{multline*}
\bigl\|(|t|-|x|+2)^\theta\, u(t,x)\bigr\|_{L^p_tL^q_x(\R\times\Omega)}\\
\leq
C\left(\|f\|_{\dot{H}_D^\gamma(\Omega)}+\|g\|_{\dot{H}_D^{\gamma-1}(\Omega)}+
\|F\|_{L^r_tL^s_x(\R\times\Omega)}+\int
\|G(t,\cdot\,)\|_{\dot{H}_D^{\gamma-1}(\Omega)}\:dt\right).
\end{multline*}
for any $\theta < 1$.
\end{lemma}

\noindent{\it Proof:} 
We will establish the result for $t\geq 0$.  
We begin this proof in the same manner as in Smith-Sogge
\cite{SS1}.  Start by observing that by \eqref{1.2} and Duhamel's
principle, the result holds for $t\in [0,1]$ and by \eqref{1.2}
\begin{multline}\label{5.1}
\|u(1,\cdot)\|_{\dot{H}_D^\gamma(\Omega)}+\|\partial_tu(1,\cdot)\|_{\dot{H}_D^{\gamma-1}(\Omega)}\\
\leq
C\left(\|f\|_{\dot{H}_D^\gamma(\Omega)}+\|g\|_{\dot{H}_D^{\gamma-1}(\Omega)}+
\|F\|_{L^r_tL^s_x(\R\times\Omega)}+\int\|G(t,\cdot)\|_{\dot{H}_D^{\gamma-1}(\Omega)}\:dt\right).
\end{multline}
By considering $t\geq 1$, finite propagation speed and support
considerations allow us to take $F=G=0$ with $f,g$ supported in
$\{|x|\leq 2\}$.

We now fix a smooth $\beta$ with $\beta(x)=1$ for $|x|\leq
\frac{1}{2}$, and $\beta(x)=0$ for $|x|\geq 1$.  We, then, write
$u$ as $u=\beta u + (1-\beta) u$ and will examine these pieces
separately.

We begin by looking at $\beta u$.  Notice that
$$\Box (\beta u)=\sum_{j=1}^n b_j(x)\partial_{x_j}u+c(x)u\equiv
\tilde{G}(t,x)$$ where $b_j,c$ are supported in $\frac{1}{2}\leq
|x|\leq 1$. Since $|t|-|x|\leq |t|$, it will suffice to show
$$\|(t+2)^{\theta}\beta(x)u(t,x)\|_{L^p_tL^q_x([1,\infty)\times\Omega)}\leq
C\left(\|f\|_{\dot{H}_D^\gamma(\Omega)}+\|g\|_{\dot{H}_D^{\gamma-1}(\Omega)}\right).$$

By the local Strichartz estimate \eqref{1.2}, Duhamel's principle,
and local energy decay \eqref{2.1}, we have
\begin{align*}
\|(t+2)^\theta&\beta(x)u(t,x)\|^p_{L^p_tL^q_x([1,\infty)\times\Omega)}\\&\leq
\sum_{j=1}^{\infty}
(j+3)^{p\theta}\|\beta(x)u(t,x)\|^p_{L^p_tL^q_x([j,j+1]\times\Omega)}\\
&\leq
C\sum_{j=1}^{\infty}(j+3)^{p\theta}\bigl(\|\beta(\cdot)u(j,\cdot)\|_{\dot{H}_D^\gamma(\Omega)}
  +\|\beta(\cdot)\partial_t
  u(j,\cdot)\|_{\dot{H}_D^{\gamma-1}(\Omega)}\\
&\quad\quad\quad\quad\quad+\int_j^{j+1}\|\tilde{G}(s,\cdot)\|_{\dot{H}_D^{\gamma-1}(\Omega)}\:ds\bigr)^p\\
&\leq C\sum_{j=1}^\infty
\frac{(j+3)^{p\theta}}{(j+3)^{p(n/2)}}\left(\|f\|_{\dot{H}_D^\gamma(\Omega)}+
  \|g\|_{\dot{H}_D^{\gamma-1}(\Omega)}\right)^p\\
&=C\left(\|f\|_{\dot{H}_D^{\gamma}(\Omega)}+\|g\|_{\dot{H}_D^{\gamma-1}(\Omega)}\right)^p
\end{align*}
so long as $\frac{n}{2}-\theta>\frac{1}{p}$.  For $n\ge 4$, since
$p\ge 2$, we have the above inequality provided $\theta <
\frac{n-1}{2}$. 

For the $v(t,x)=(1-\beta)(x)u(t,x)$ piece, we have that $v$
satisfies the Minkowski wave equation
$$\begin{cases}
\Box v(t,x)=-\tilde{G}(t,x)\\
v(0,x)=(1-\beta)(x)f(x)\\
\partial_tv(0,x)=(1-\beta)(x)g(x).\end{cases}$$
\par

Write $v=v_0+v_1$ where $v_0$ solves the homogeneous wave equation
with the same Cauchy data as $v$ and $v_1$ solves the
inhomogeneous wave equation with vanishing Cauchy data.  Then, by
Theorem \ref{theorem3.4}, we have
\begin{multline*}
\|(t-|x|+2)^{\theta}(1-\beta)(x)u(t,x)\|_{L^p_tL^q_x(\R\times\Omega)}\\
\le
C\left(\|(1-\beta)(x)f\|_{\dot{H}^{\gamma}(\R^n)}+\|(1-\beta)(x)g(x)\|_{\dot{H}^{\gamma-1}(\R^n)}\right)
  \\+\|(t-|x|+2)^{\theta} v_1(t,x)\|_{L^p_tL^q_x(\R\times\Omega)}.
\end{multline*}

When $n\ge 4$, we can handle the last piece easily using local
energy decay.  By Duhamel's principle, write
$$\|(t-|x|+2)^{\theta}
v_1(t,x)\|_{L^p_tL^q_x(\R\times\Omega)}=\Bigl\|(t-|x|+2)^{\theta}\int_0^{t}
v_1(s;t-s,x)\:ds\Bigr\|_{L^p_tL^q_x(\R\times\Omega)}$$ where
$v_1(s;\cdot,\cdot)$ solves
$$\begin{cases}
\Box v_1(s;t,x)=0\\
v_1(s;0,x)=0\\
\partial_tv_1(s;0,x)=-\tilde{G}(s,x).\end{cases}$$
Applying Minkowski's integral inequality, we have
$$\|(t-|x|+2)^{\theta} v_1(t,x)\|_{L^p_tL^q_x(\R\times\Omega)}
\le C\int s^{\theta} \|(t-s-|x|+2)^{\theta}
v_1(s;t-s,x)\|_{L^p_tL^q_x(\R\times\Omega)}\:ds.$$ Thus, by
Theorem \ref{theorem3.4}, the right side is bounded by
$$\int s^{\theta} \|\tilde{G}(s,\cdot)\|_{\dot{H}^{\gamma-1}(\R^n)}\:ds.$$
Finally, by local energy decay \eqref{2.1}, we have that this is bounded
by
$$C\int s^{\theta-\frac{n}{2}}
\left(\|f\|_{\dot{H}^{\gamma}(\Omega)}+\|g\|_{\dot{H}_D^{\gamma-1}(\Omega)}\right)\:ds\le
C\left(\|f\|_{\dot{H}^{\gamma}(\Omega)}+\|g\|_{\dot{H}_D^{\gamma-1}(\Omega)}\right)$$
for $\theta<1$.

\qed

We are now ready to prove the main theorem.

\noindent{\it Proof of Theorem \ref{theorem1.1}}. By the previous
lemma, it will suffice to show the result when $f$ and $g$ vanish
for $\{|x|\le 1\}$.
\par
We begin by decomposing $u$ into
$$u(t,x)=u_0(t,x)-v(t,x)$$
where $u_0$ solves the Minkowski wave equation
$$\begin{cases}
\Box u_0(t,x)=F(t,x)\\
u_0(0,x)=f(x)\\
\partial_tu_0(0,x)=g(x).
\end{cases}$$
Here $F$ is assumed to be $0$ on $\R^n\backslash\Omega$.
\par
We now fix a smooth compactly supported $\beta$ such that
$\beta(x)=1$ for $|x|\le 1/2$ and $\beta(x)=0$ for $x\ge 1$. Then,
further decompose $u$ into
$$u(t,x)=u_0(t,x)-v(t,x)=(1-\beta)(x)u_0(t,x)+(\beta(x)u_0(t,x) - v(t,x)).$$
By the Global Minkowski Strichartz estimate \eqref{1.3},
$(1-\beta)(x)u_0(t,x)$ satisfies the desired estimate.  Thus, we
may focus on $\beta(x)u_0(t,x)-v(t,x)$.
\par
We have that $\beta(x)u_0(t,x)-v(t,x)$ satisfies
$$\Box (\beta(x)u_0(t,x)-v(t,x))=\beta(x)F(t,x)+G(t,x)$$
with zero Cauchy data (since we are assuming that $f, g$ vanish
for $|x|\le 1$).  Here
$$G(t,x)=\sum_{j=1}^n b_j(x)\partial_{x_j}u_0(t,x)+c(x)u_0(t,x).$$
where $b_j, c$ vanish for $|x|\ge 1$.  By Lemma \ref{lemma4.2},

\begin{equation}\label{5.3}
\int_{-\infty}^\infty
\|G(t,\cdot)\|^2_{\dot{H}^{\gamma-1}_D(\Omega)}\le
C\left(\|f\|_{\dot{H}^{\gamma}_D(\Omega)} +
\|g\|_{\dot{H}^{\gamma-1}_D(\Omega)}+\|F\|_{L^r_tL^s_x(\R\times\Omega)}\right)^2\end{equation}

\par
Let
\begin{align*}
F_j(t,x)&=\chi_{[j,j+1]}(t)F(t,x)\\
G_j(t,x)&=\chi_{[j,j+1]}(t)G(t,x),
\end{align*}
and write (for $t>0$)
$$\beta u_0-v = \sum_{j=0}^{\infty} u_j(t,x)$$
where $u_j(t,x)$ is the forward solution to
$$\Box u_j(t,x)=\beta(x)F_j(t,x)+G_j(t,x)$$
with zero Cauchy data.
\par
Thus, by Lemma \ref{lemma5.1}, we have
\begin{multline}\label{5.4}
\|(t-j-|x|+2)^{\theta}u_j(t,x)\|_{L^p_tL^q_x(\R\times\Omega)}\\\le
  C\left(\|F_j(t,x)\|_{L^r_tL^s_x(\R\times\Omega)}
  +\int_{j}^{j+1}\|G(t,\cdot)\|_{\dot{H}^{\gamma-1}_D(\Omega)}\:dt\right).\end{multline}
\par
Since $u_j$ is supported in the region $t-j-|x|+2\ge 1$, an
application of the Cauchy-Schwartz inequality yields
\begin{align*}
|\beta(x)&u_0(t,x)-v(t,x)|\le \sum_{j=0}^{\infty} |u_j(t,x)|\\
&\le \left(\sum_{j=0}^{\infty}
  (t-j-|x|+2)^{-2\theta}\right)^{1/2}\:\left(\sum_{j=0}^{\infty}
  [(t-j-|x|+2)^{\theta}u_j(t,x)]^2\right)^{1/2}\\
&\le C\left(\sum_{j=0}^{\infty}
  [(t-j-|x|+2)^{\theta}u_j(t,x)]^2\right)^{1/2}
\end{align*}
since we can choose $\theta>1/2$.
\par
Since $1\le r,s\le 2 \le p,q$, Minkowski's integral inequality,
\eqref{5.3}, and \eqref{5.4} yield
\begin{align*}
\|\beta(x)u_0(t,x)&-v(t,x)\|^2_{L^p_tL^q_x(\R\times\Omega)}\\
&\le
C\sum_{j=0}^{\infty}\|(t-j-|x|+2)^{\theta}u_j\|^2_{L^p_tL^q_x(\R\times\Omega)}\\
&\le
C\sum_{j=0}^{\infty}\|F_j\|^2_{L^r_tL^s_x(\R\times\Omega)}+C\sum_{j=0}^\infty
\left(\int_j^{j+1}\|G(t,\cdot)\|_{\dot{H}_D^{\gamma-1}(\Omega)}\:dt\right)^2\\
&\le
C\sum_{j=0}^{\infty}\|F_j\|^2_{L^r_tL^s_x(\R\times\Omega)}+C\sum_{j=0}^\infty
\left(\int_j^{j+1}\|G(t,\cdot)\|^2_{\dot{H}_D^{\gamma-1}(\Omega)}\:dt\right)\\
&\le C\|F\|^2_{L^r_tL^s_x(\R\times\Omega)}+C\int_0^{\infty}\|G(t,\cdot)\|^2_{\dot{H}_D^{\gamma-1}(\Omega)}\:dt\\
&\le
C\left(\|f\|_{\dot{H}^{\gamma}_D(\Omega)}+\|g\|_{\dot{H}_D^{\gamma-1}(\Omega)}+\|F\|_{L^r_tL^s_x(\R\times\Omega)}\right)^2
\end{align*}
as desired. \qed


\end{document}